\documentclass[12pt]{amsart}

\usepackage{amsmath}
\usepackage{amssymb,mathrsfs}
\usepackage[left=3.5cm,right=3.5cm,top=3.5cm,bottom=3.5cm]{geometry}
\usepackage{graphicx}

\newcommand{\eps}{\varepsilon}
\newcommand{\al}{\alpha}
\newcommand{\be}{\beta}
\newcommand{\ga}{\gamma}

\begin{document}

\title[The method and examples]{The method and examples of solving problems analogous to
``the problem of three bisectors''}

\author{S. F. Osinkin }
\address{ Moscow State University, Russia}
\email[S. F. Osinkin]{osinkin1947@yandex.ru}

\begin{abstract}
We suggest a method of solving the problem of existence of a triangle with prescribed 
two bisectors and one third element which can be taken as one of the angles, 
the sides, the heights or the medians, or the third bisector.
\end{abstract}


\maketitle

{\em In memory of our school mathematics teacher Maya Aronovna Kaganovskaya-Sazonova}

\section*{Introduction}

In the Geometry of Triangles one can distinguish a group of problems with similar
formulations which have not been discussed in popular mathematics literature with
attention that they deserve. Apparently, this indifference was caused by
lack of a simple enough method of solving these problems.
As an exception one can consider ``the problem of three bisectors'' which was
formulated in 1875 by French mathematician H.~Brocard (1845-1922) in the form of questions:

{\it (i) Does there exist a triangle with prescribed lengths $l_a,l_b,l_c$ of its bisectors?

(ii) Is it necessary to impose any conditions on these lengths for the existence of such a triangle?}

Solution of this problem, based on the Brower fixed point theorem, was given by Romanian
mathematicians  P.~Mironescu and L.~Panaitopol in Ref.~\cite{mp-1994}. Later two more other solutions
of this problem were suggested in Refs.~\cite{za-2003,os-2016}. However, if we replace in the Brocard
formulation of the problem one of bisectors by another element of a triangle, say, by its side, or height,
or median, {\it etc.}, then we arrive at the necessity of finding the general approach to this type of
problems. Such an approach does exist and the aim of this article is to demonstrate by solution of
several concrete problem the efficiency of the method as well as its simplicity and visual clearness.

Let us discuss briefly the essence of the method. Let we consider the question of existence of a triangle
$ABC$ with three prescribed elements, for example, its median $m_a$, height $h_b$, and bisector $l_c$.
We choose from these three given parameters any two, assuming that they satisfy the necessary conditions
(for example, $m_a\geq h_b/2$). Then we introduce into consideration some angle in such a way, that
it would be possible to construct a triangle with our two chosen elements and this angle (say, in our example
we have $m_a$, $h_b$, and $\angle C$). Thus, we parameterize triangles by the values of these three elements
which can be considered as variables of the problem. As we shall see below, depending on the problem under
consideration, we can fix one or two of these parameters. For example, in particular case of the problem with
given heights $h_a$ and $h_b$ ($h_a\leq h_b$) the parameter $h_a$ can be considered as a variable changing
from $h_a=0$ to $h_a=h_b$, but the value of $h_b$ is fixed ($h_b=1$). However, in the problem with given
two bisectors $l_a$ and $l_b$ they both can be treated as constants. In any case, if we fix one prescribed
parameter and fix the value of the parameter considered as a variable, then we can calculate the other elements
as functions of the angle introduced into consideration. Changing the angle in the admissible interval of
its values, we determine the interval of change of the third element as a function of two or one varying elements.
As a result, we find the relationship between all three elements which formulates the necessary and
sufficient conditions for existence of a triangle under consideration.

To determine, how many triangles with given values of three elements exist, we should study if the
function which defines the third element is a continuous function of the angle introduced into our
solution. If it is monotonous, then the triangle is unique. If it is not, then the number
of such triangles depends on the value of the third element. For example, in case of the elements
$h_a, h_b$ and $l_a$ ($h_a\leq h_b$) we can get one, two or three solutions depending on the values
of $h_a$ and $l_a$.

\section{Triangles with two prescribed heights}

\begin{figure}[ht]
\begin{center}
\includegraphics[width=8cm]{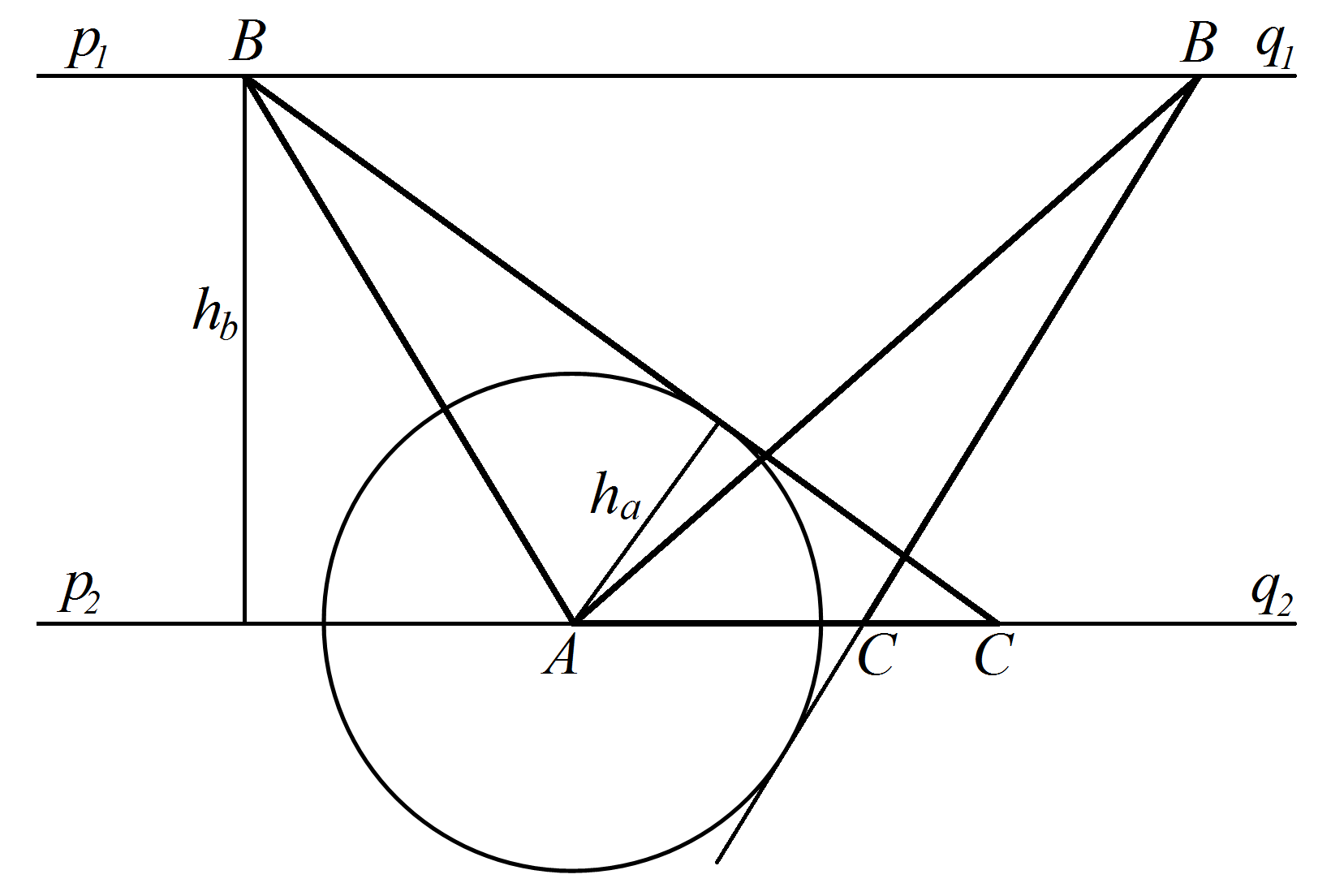}
\caption{Construction of triangles $ABC$ with prescribed  lengths of two heights $h_a=h_1$ and
$h_b=h_2$ ($h_1\leq h_2$). }
\end{center}\label{fig1}
\end{figure}

\begin{figure}[ht]
\begin{center}
\includegraphics[width=6cm]{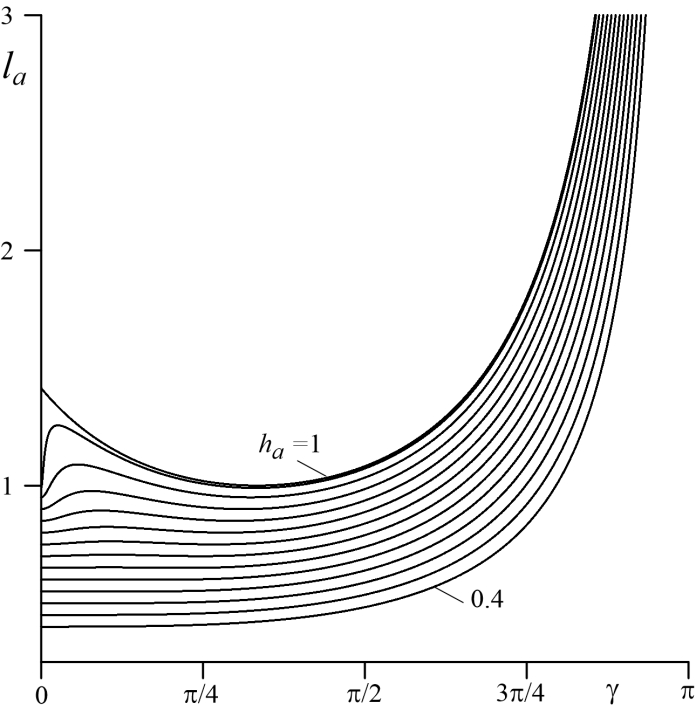}
\put(-100,160){(a)}
\hspace{5mm}
\includegraphics[width=6cm]{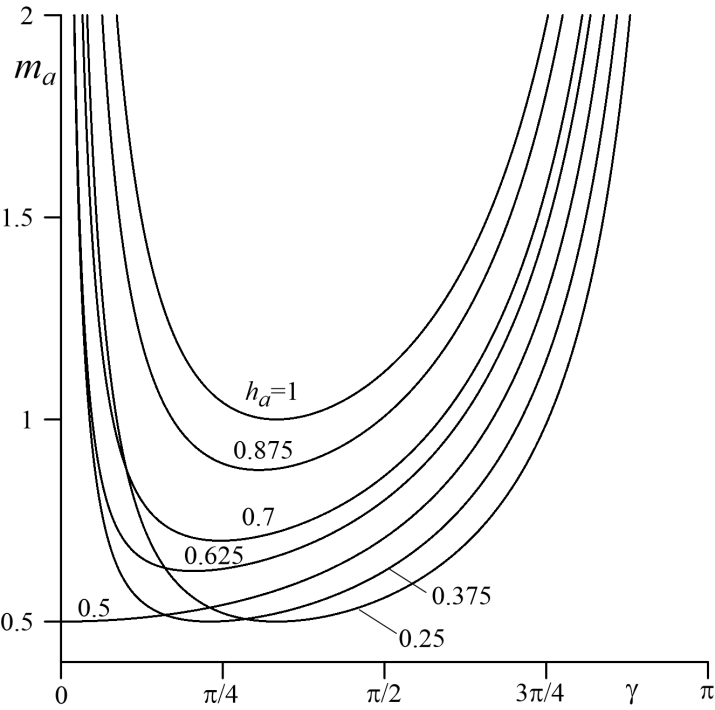}
\put(-100,160){(b)}
\caption{Dependence of the bisector's length $l_q$ (a) and of the median's length $m_a$
on the angle $\gamma$ for different values of $h_a$ ($0< h_a\leq h_b$) and $h_b=1$.
}
\label{fig2}
\end{center}
\end{figure}

Let us consider a set of triangles $ABC$ (for convenience we use the same notation for vertices of
all triangles) with prescribed lengths of heights $h_a=h_1$ and $h_b=h_2$ with $h_1\leq h_2$.
This set can be constructed in the following way: (i) we draw two parallel lines $p_1q_1$ and
$p_2q_2$ at distance $h_2$ between them (see Fig.~\ref{fig1}); (ii) we take the point $A$ in the lower
line $p_2q_2$ (which will be the common vertex of all triangles under consideration) and draw the
circle with the radius $h_1$ and with the center at this point; (iii) at last, we draw all possible lines
which touch the circle at the points in its right part, as is shown in Fig.~\ref{fig1}.
We denote the intersection points of the tangent lines with the lines $p_1q_1$ by the symbols $B$ and $C$
(there are infinitely many of such pairs of the points) and connect $A$ with $B$ and $B$ with $C$, so that
we get the triangles $ABC$ with $h_a=h_1$ and $h_b=h_2$. In this construction, the vertex $C$ is always
located  to the right of the vertex $A$ independently of the tangency point. Now we introduce the
angle $\gamma=\angle C$. It is evident that when we change the location of the tangency point from the
lowest position at the vertical diameter edge to its highest position, the angle $\gamma$ varies
continuously from $\gamma=\pi$ to $\gamma=0$. If the sides of the triangle $ABC$ are denoted as $a,b$ and $c$,
then we get
$$
a=\frac{h_b}{\sin\gamma},\qquad b=\frac{h_a}{\sin\gamma},\qquad c=\sqrt{a^2+b^2-2ab\cos\gamma}.
$$
Thus, if we know $h_a$, $h_b$ and $\gamma$, then it is easy to calculate the lengths of the sides of
the triangle $ABC$ and, consequently, the other its elements --- heights, medians, bisectors, {\it etc.}
It is especially interesting to study the dependence of the bisector $l_a$ and the median $m_a$ on
$h_a$ and $h_b$. We have given in Fig.~\ref{fig2} the plots of the functions $l_a(\gamma)$ (a) and
$m_a(\gamma)$ (b) for $h_2=1$ and some set of values of $h_1$.

For this plotting we have used the formulas
\begin{equation}\nonumber
  \begin{split}
  & l_a=\left\{bc\left[1-\frac{a^2}{(b+c)^2}\right]\right\}^{1/2},\\
  & m_a=\frac12\left(2b^2+2c^2-a^2\right)^{1/2}=\frac12\left(4b^2+a^2-4ab\cos\gamma\right)^{1/2}.
  \end{split}
\end{equation}
From these plots and formulas one can find that $l_a(\gamma)$ takes the minimal value
$l_a=h_a=h_1$ at two values of $\gamma_1=0$ and $\gamma_2=\arccos({h_2}/({2h_1}))$ or
at only one value $\gamma_1=0$ or $\gamma_2=\pi/3$, depending on the value of $h_1$:
if $\tfrac12h_2<h_1<h_2$, then we get the case with two minima, if $h_1\leq \tfrac12h_2$
or $h_1=h_2$, then we have a single minimum. The value $\gamma_2$ corresponds to an
isosceles triangle with $b=c$ (in this case the tangent point of the side $BC$ with the circle
in Fig.~\ref{fig1} is located at the same distances from the lines $p_1q_1$ and $p_2q_2$, so
consequently $m_a=h_a$ and $b=c$). The value $\gamma_1$ corresponds to a triangle with infinitely
large sides. Therefore an actual triangle does exist for however small value of
$\gamma$ but not for $\gamma=0$, although the bisector $l_a(0)=h_1$ in this case has
the prescribed value, too.
It is obvious that with decrease of $h_1$ starting from $h_1=h_2$, the angle $\gamma_2$
decreases also from the value $\gamma_2=\tfrac13\pi$ reaching $\gamma_2=0$ at $h_1=\tfrac12h_2$.
Hence for $h_1\leq \tfrac12h_2$ the minimum of $l_a$ corresponding to an isosceles triangle
stops to exist and the dependence $l_a(\gamma)$ can be represented in the form
$$
l_a(\gamma)=h_1\left(1+\frac{k\gamma^2}{\pi-\gamma}\right),
$$
where $0\leq\gamma<\pi$ and $k$ is some constant parameter.

It is clearly seen in Fig.~\ref{fig2}(a) that the function $l_a(\gamma)$ has a local maximum
in the interval between $\gamma_1$ and $\gamma_2$. Obviously, this maximum exists as long as
there are two local minima of $l_a(\gamma)$ at $\gamma_1$ and $\gamma_2$, that is for
$\tfrac12h_2<h_1<h_2$. This maximum value of $l_a$ depends on $h_a$ --- the smaller is the
difference between $h_b$ and $h_a$ ($h_a\leq h_b$), to larger extent a part of the curve
$l_a(\gamma)$ in vicinity of this maximum enters into the angle formed by the ordinate axis and
the curve $l_a(\gamma)$ drawn for $h_a=h_b$. Thus, with increase of $h_a$ from $h_a=\tfrac12h_b$
to $h_a=h_b$, the value of the local maximum $l_{\text{max}}(h_a)$ increases monotonously from
$\tfrac12h_b$ to $\sqrt{2}\,h_b$. The edge points of the interval within which $l_{\text{max}}(h_a)$
changes (i.e., $\tfrac12h_b$ and $h_b$) do not correspond to any real triangles, because the
corresponding values of $\gamma$ are equal to zero. Consequently, for any $h_a$ in the interval
$(\tfrac12h_b,h_b)$ and any $l_a$ from the interval $(h_a,l_{\text{max}}(h_a))$ the straight line
$l_a(\gamma)=l_a$ parallel to the abscissa axis intersects the graph of the function $l_a(\gamma)$
in Fig.~\ref{fig2}(a) in three points with different values of $\gamma$. This meand that for these
$h_a$ and $l_a$ there exist three triangles with the same values of $h_a$, $h_d$ and $l_a$.
Let us prove that these triangles are not equal to each other. We suppose the opposite, say, that
$A_1B_1C_1=A_2B_2C_2$. Since $h_{a1}=h_{a2}=h_1$, $h_{b1}=h_{b2}=h_2$ and in equal triangle the
equal heights correspond to equal sides, we have $B_1C_1=B_2C_2$ and $A_1C_1=A_2C_2$. From
equality of triangles it follows that also $A_1B_1=A_2B_2$, and then $\angle C_1=\angle C_2$,
what contradicts to our condition.

Besides this situation, two other are possible:

(a) there are two solutions (triangles), if $h_a=h_b$ and $h_b< l_a<\sqrt{2}\, h_b$
or $\tfrac12h_b< h_a< h_b$ and $l_a=l_{\text{max}}(h_a)$;

(b) there is only one solution (triangle), if (i) $h_a=h_b$ and $l_a\geq \sqrt{2}\,h_b$,
(ii) $\tfrac12h_b<h_a<h_b$ and $l_a>l_{\text{max}}(h_a)$ or $l_a=h_a$;
(iii) $h_a\leq \tfrac12h_b$ and $l_a>h_a$.

To complete the discussion, we provide the formulas for dependence of $l_a(\gamma)$ in
vicinity of $\gamma=0$:
\begin{equation}\nonumber
  \begin{split}
  & l_a(\gamma) =\sqrt{2}\,h_2\left(1-\frac14\gamma\right),\qquad \text{if}\qquad h_1=h_2,\\
  & l_a(\gamma)=h_1\left[1+\frac{2h_1^2+2h_2^2-h_1h_2}{12(h_2-h_1)^2}\cdot\gamma^2\right],
  \qquad\text{if}\qquad h_1\neq h_2.
  \end{split}
\end{equation}

Besides that, the straight line $\gamma=\pi$ is a vertical asymptote for all curves $l_a(\gamma)$
independently of the value of $h_a$.

Thus, if there are three segments which satisfy the conditions $S_1\leq S_2\leq S_3$, then there
always exist one, two or three triangles with $h_a=S_1$, $h_b=S_2\,(S_3)$ and $l_a=S_3\,(S_2)$.

Now we turn  to another case. We show in Fig.~\ref{fig2}(b) the plots of the function $m_a(\gamma)$,
$$
m_a(\gamma)=\frac1{2\sin\gamma}(4h_1^2+h_2^2-4h_1h_2\cos\gamma)^{1/2},
$$
for values of the parameter $h_1$ from $h_1=h_2=1$ to $h_1=0.25$. It is clear that for each $h_1$ the
function $m_a(\gamma)$ has a single minimum which depends on $h_1$, if $\tfrac12h_2\leq h_1\leq h_2$,
and does not depend on $h_1$ for $h_1\leq\tfrac12h_2$. Besides that, if $h_1\neq\tfrac12h_2$,
then each curve $m_a(\gamma)$ has two vertical asymptotes $\gamma=0$ and $\gamma=\pi$. But in
the case $h_1=\tfrac12h_2$ the asymptote $\gamma=0$ disappears and there remains the asymptote $\gamma=\pi$ 
only. Notice that all curves in Fig.~\ref{fig2}(b) can be subdivided into two families: the first family
consists of the curves which do not cross each other and correspond to $h_1$ from the interval
$[\tfrac12h_2,\,h_2]$, and the second family contains the curves which cross each other and the
curves of the first family and they correspond to $h_1<\tfrac12h_2$. Each point of crossing of
two curves in Fig.~\ref{fig2}(b) means that there exist two different triangles with equal values of
$m_a$, $h_b$ and the angle $\gamma$. The absence of such crossing poins means that there exists
only one such a triangle. For example, the curves in Fig.~\ref{fig2}(a) do not cross each other and hence
there exists only one triangle with given values of $l_a$, $h_b$ and the angle $\gamma$.

Let us find the minimal value of $m_a(\gamma)$. To this end, we calculate the derivative
$$
m_a'(\gamma)=\frac1{4m_a\sin^3\gamma}[2h_1h_2\cos^2\gamma-(4h_1^2+h_2^2)\cos\gamma+2h_1h_2]
$$
and find roots of the equation $m_a'(\gamma)=0$. After simple transformations we find that for any
$h_1$ from the interval $[\tfrac12h_2,h_2]$ the minimum is given by
$\mathrm{min}\,m_a(\gamma)=h_1$ at $\gamma=\gamma_m=\arccos(h_2/(2h_1))$, but if $h_1\leq\tfrac12h_2$,
then $\mathrm{min}\,m_a(\gamma)=\tfrac12h_2$ and this value does not depend on $h_1$ at the points
of minima with $\gamma_m=\arccos(2h_1/h_2)$. These results show that the angle $\gamma_m$ decreases
from $\pi/3$ to zero when we decrease $h_1$ from $h_1=h_2$ to $h_1=\tfrac12h_2$ and it increases
from zero to $\pi/2$ when we decrease $h_1$ further from $h_1=\tfrac12h_2$ to arbitrary small value.
If $h_1=\tfrac12h_2$, then in vicinity of $\gamma=0$ the function $m_a(\gamma)$ can be
approximated by
$$
m_a(\gamma)\approx \frac12h_2\left(1+\frac18\gamma^2\right).
$$

The above analysis shows that there exist one or two triangles with $h_a=S_1$, $h_b=S_2$ and
$m_a=S_3$, where $S_1,S_2,S_3$ are given segments which lengths satisfy the inequalities
$S_1\leq S_2\leq S_3$.

\section{On existence of triangles with prescribed lengths of height, median and bisector
which are referred to the different angles}

\begin{figure}[ht]
\begin{center}
\includegraphics[width=9cm]{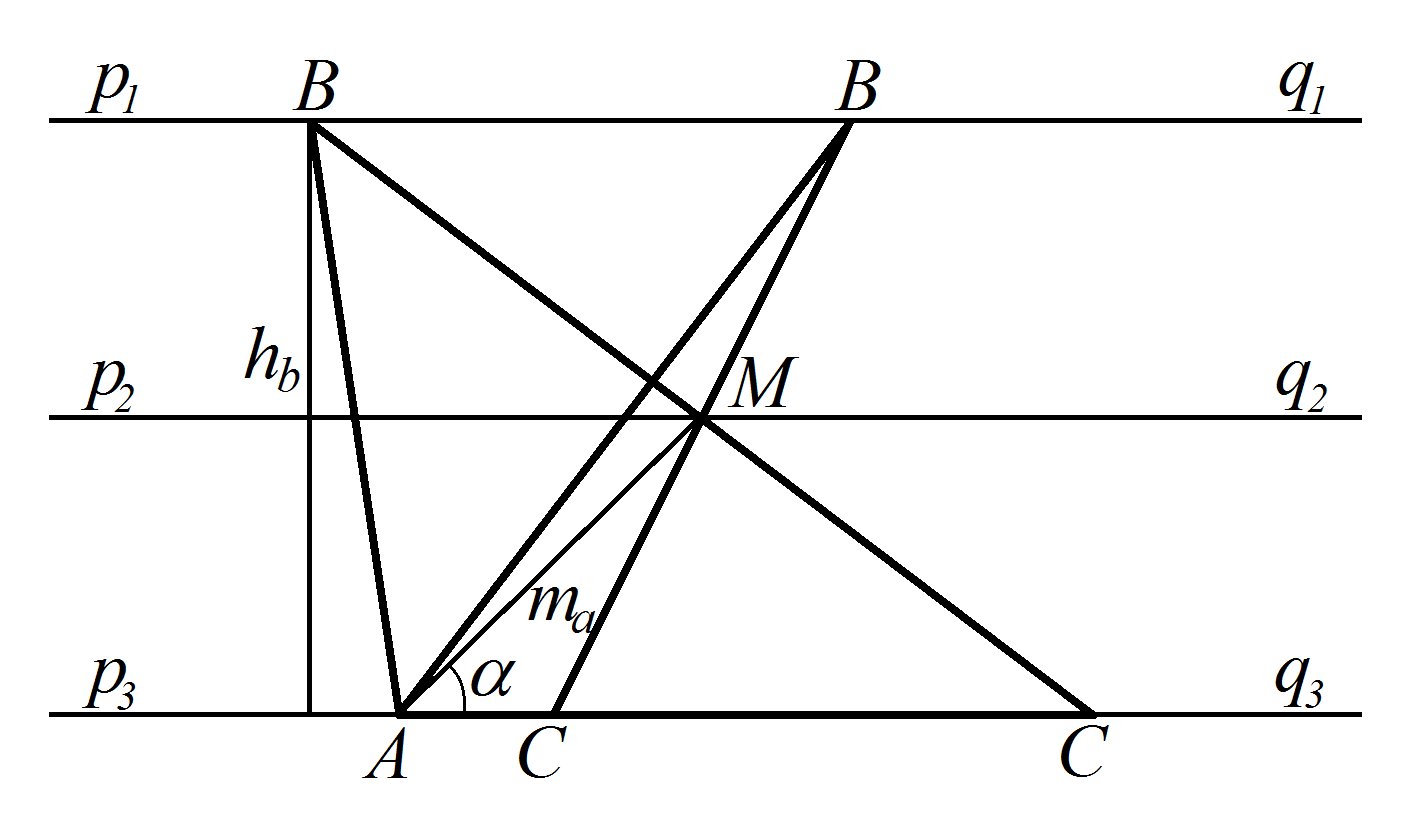}
\caption{Construction of triangles with prescribed lengths of the median $m_a$ and the height $h_b$.
}
\label{fig3}
\end{center}
\end{figure}

Let us take three arbitrary segments with lengths $S_1,S_2,S_3$ and let we have $\tfrac12S_1\leq S_2$.
We draw three parallel lines $p_1q_1$, $p_2q_2$ and $p_3q_3$ with distance between them equal to
$\tfrac12S_1$, as is shown in Fig.~\ref{fig3}. We take an arbitrary point $A$ in the line $p_3q_3$
and join it by the segment $S_2$ with the point $M$ lying in the line $p_2q_2$. This always can be
done provided $\tfrac12S_1\leq S_2$. We denote the slope of the line $S_2$ with respect to the line
$p_3q_3$ by the angle $\alpha$ directed counterclockwise. The angle $\alpha$ can have two
essentially different values (excluding the value $\alpha=\pi/2$): $\alpha=\alpha_1<\tfrac12\pi$
and $\alpha=\alpha_2=\pi-\alpha_1>\tfrac12\pi$. We draw through the point $M$ the line which
intersects $p_1q_1$ and $p_3q_3$ at the points $B$ and $C$, respectively. Then, we join $A$ and $B$
by the segment and obtain the triangle $ABC$ with $m_a=AM=S_2$ and $h_b=S_1$. If we have a set
of lines drawn through the point $M$ such that the angle $\angle AMC$ increases continuously
from zero to $\pi-\alpha$, then we obtain a set of triangles with equal lengths of $m_a$ and $h_b$.
Let us denote the set of triangles with $\alpha=\alpha_1$ as $\aleph_1$, the set of triangles with
$\alpha=\alpha_2$ as $\aleph_2$, and the set of triangles with $\alpha=\tfrac12\pi$ as $\aleph$.

Now we show that when we change the angle $\angle AMC$ from zero to $\pi-\alpha$, then the
bisector $l_c$ of the angle $\angle C$ increases monotonously from $l_c=0$ to $l_c=\infty$.
It is evident that for $\alpha=\alpha_2$ and $\alpha=\tfrac12\pi$ the sides $AB$ and $BC$
increase monotonously up to infinity with growth of the angle $\angle AMC$, whereas the
angle $\angle C$ decreases monotonously down to zero. We find from the corresponding formulas
$$
l_c=\frac{2AC\cdot BC}{AC+BC}\cos\tfrac12\angle C=\frac{2\cos\tfrac12\angle C}{\frac1{AC}+\frac1{BC}}
\qquad\text{and}\qquad
l_c=2AC\cdot\cos\tfrac12\angle C,
$$
that for however small $AC$ the length $l_c$ increases up to infinity. Besides that, since at $\angle AMC=0$
we have $AC=0$, we conclude that $l_c$ increases monotonously from $l_c=0$ to $l_c=\infty$.

If $\alpha=\alpha_1$,
then with growth of $\angle AMC$ the length of $DC$ at first decreases from $BC=2AM$ to $BC=h_b$ and
after that increases up to infinity. To prove that $l_c$ increases monotonously, we have to find the
derivative of $l_c$ with respect to the varying angle $\angle C$ and to check that it is always positive.
Let us denote $\angle C=\gamma$, $h_b=h$, $m_a=m$, so that $\sin\alpha=\frac{h}{2m}$ and
$BC=\frac{h}{\sin\gamma}$. From $AC/\sin(\alpha+\gamma)=m/\sin\gamma$ we find
$AC=m\sin(\alpha+\gamma)/\sin\gamma$ and, after substitution of $AC$ and $BC$ into the
formula for $l_c$, we get
\begin{equation}\nonumber
  \begin{split}
  l_c&=\frac{2mh\sin(\alpha+\gamma)\cos\tfrac12\gamma}{\sin\gamma\,[h+m\sin(\alpha+\gamma)]}=
  \frac{mh\sin(\alpha+\gamma)}{\sin\left(\tfrac12\gamma\right)[h+m\sin(\alpha+\gamma)]}\\
  &=\frac{mh\sin(\alpha+\gamma)}{m\sin\left(\tfrac12\gamma\right)[\frac{h}m+\sin(\alpha+\gamma)]}=
  \frac{h\sin(\alpha+\gamma)}{\sin\left(\tfrac12\gamma\right)[2\sin\alpha+\sin(\alpha+\gamma)]}.
  \end{split}
\end{equation}
On the contrary to the behavior of the angle $\angle AMC$, which changes from zero to $\pi-\alpha$, 
the angle $\gamma$ changes from $\pi-\alpha$ to zero. Let us calculate the derivative $l'_c(\gamma)$
and show that $l'_c(\gamma)<0$. For convenience, we calculate the logarithmic derivative
equal to
$$
l'_c(\gamma)=l_c(\gamma)[\ln l_c(\gamma)]'.
$$
We introduce the notation $f=[\ln l_c(\gamma)]'$ and find
\begin{equation}\nonumber
\begin{split}
  f&=\left\{\ln\left[\frac{h\sin(\alpha+\gamma)}
  {\sin\frac{\gamma}2(2\sin\alpha+\sin(\alpha+\gamma))}\right]\right\}'\\
  &=\frac{\cos(\alpha+\gamma)}{\sin(\alpha+\gamma)}-\frac12\frac{\cos\frac{\gamma}2}{\sin\frac{\gamma}2}-
  \frac{\cos(\alpha+\gamma)}{2\sin\alpha+\sin(\alpha+\gamma)}.
  \end{split}
\end{equation}
This trigonometric expression can be cast to the form
\begin{equation}\nonumber
  f=\frac{2\sin\alpha\left[\sin\frac{\gamma}2\cdot\cos(\alpha+\gamma)
  -\sin\left(\alpha+\frac{\gamma}2\right)\right]-\sin^2(\alpha+\gamma)\cos\frac{\gamma}2}
  {2\sin\frac{\gamma}2\cdot\sin(\alpha+\gamma)[2\sin\alpha+\sin(\alpha+\gamma)]}.
\end{equation}
We denote
\begin{equation}\nonumber
  \begin{split}
  f_1 & =\frac{2\sin\alpha\left[\sin\frac{\gamma}2\cdot\cos(\alpha+\gamma)
  -\sin\left(\alpha+\frac{\gamma}2\right)\right]}
  {2\sin\frac{\gamma}2\cdot\sin(\alpha+\gamma)[2\sin\alpha+\sin(\alpha+\gamma)]}=\frac{u}v,\\
  f_2 & =\frac{-\sin^2(\alpha+\gamma)\cos\frac{\gamma}2}
  {2\sin\frac{\gamma}2\cdot\sin(\alpha+\gamma)[2\sin\alpha+\sin(\alpha+\gamma)]}.
  \end{split}
\end{equation}
Since $\sin\frac{\gamma}2\geq0$ and $\cos(\alpha+\gamma)\leq1$, then
$$
u\leq 2\sin\alpha\left[\sin\frac{\gamma}2-\sin\left(\alpha+\frac{\gamma}2\right)\right]=
-4\sin\alpha\sin\frac{\alpha}2\cos\frac{\alpha+\gamma}2.
$$
Consequently,
\begin{equation}\nonumber
\begin{split}
\frac{u}v & \leq\frac{-4\sin\alpha\sin\frac{\alpha}2\cos\frac{\alpha+\gamma}2}
{\sin\frac{\gamma}2\cdot\sin(\alpha+\gamma)[2\sin\alpha+\sin(\alpha+\gamma)]}\\
& = \frac{-2\sin\alpha\sin\frac{\alpha}2}
{\sin\frac{\gamma}2\cdot\sin\frac{\alpha+\gamma}2[2\sin\alpha+\sin(\alpha+\gamma)]}
=\frac{u_1}{v_1}.
\end{split}
\end{equation}
It is evident that $u_1<0$ and $v_1>0$, therefore $f_1<0$ and so by virtue of $f_2\leq0$
we get $f_1+f_2<0$
and $l'_c(\gamma)<0$. Therefore with decrease of $\gamma$ the bisector $l_c$
increases monotonously. From monotonous change of $l_c(\gamma)$ from $l_c=0$ to
$l_c=\infty$ it follows that if $\alpha\neq\pi/2$, then in each set $\aleph_1$ and
$\aleph_2$ there exists one triangle with given lengths $m_a$, $h_b$ and $l_c$
($m_a\geq \frac12h_b$). In the case $\alpha=\frac{\pi}2$ ($m_a=h_b/2$) such a triangle
exists in $\aleph$.

Let us prove that neither of triangles which belong to $\aleph_1$ does not coincide with any
triangle from $\aleph_2$ with exception of the case when $A_1C_1=A_2C_2=0$, so that both
triangles degenerate into a segment with the length $2m$. The proof is achieved by reducing
to contradiction. We assume that there exist two equal triangles $A_1B_1C_1$ from $\aleph_1$
and $A_2B_2C_2$ from $\aleph_2$. Since $h_{b1}=h_{b2}=h$ and the triangles $A_1B_1C_1$ and
$A_2B_2C_2$ are equal to each other, then $A_1C_1=A_2C_2(\neq0)$. There are two possibilities
to superimpose one triangle over the other: (i) apex $A_1$ coincides with $A_2$, $B_1$ with $B_2$,
and $C_1$ with $C_2$, or (ii) $A_1$ coincides with $C_2$, $B_1$ with $B_2$, and $C_1$ with $A_2$.
But in the first case the medians drawn from the angles $\angle A_1$ and $\angle A_2$ cannot
superimpose because of different slopes with respect to the bases of the corresponding triangles.
For realization of the second case it is necessary to have $m_{a1}=m_{c2}=m$, but since
$m_{a1}=m_{a2}$, we have also $m_{a2}=m_{c2}=m$. Consequently, the triangle $A_2B_2C_2$
is isosceles with equal angles $\angle A_2$ and $\angle C_2$, what is impossible because
$\angle A_2>\alpha_2>\pi/2$. Thus, our supposition that two triangles from $\aleph_1$ and
$\aleph_2$ are equal to each other is wrong.

Our analysis leads to the following conclusions:

1) If three segments $S_1,S_2,S_3$ are given and two of them, say, $S_1$ and $S_2$ satisfy the
inequality $\tfrac12S_1<S_2$, then there exist two different triangles $A_1B_1C_1$ and $A_2B_2C_2$
with $h_{b1}=h_{b2}=S_1$, $m_{a1}=m_{a2}=S_2$ and $l_{c_1}=l_{c2}=S_3$, and in these triangles
the median $m_{a1}$ has a slope angle $\alpha_1=\alpha$ with respect to the base $A_1C_1$, whereas
the median $m_{a2}$ has a slope angle $\alpha_2=\pi-\alpha$ with respect to $A_2C_2$, where
$\alpha=\arcsin(S_1/(2S_2))$ and $0<\alpha<\pi/2$.

2) If $\tfrac12S_1=S_2$ and, consequently, $\alpha=\pi/2$, then there exists a single triangle
$ABC$ with $h_b=S_1$, $m_a=S_2$ and $l_c=S_3$.

\section{Triangles with prescribed lengths of two bisectors}

In Ref.~\cite{os-2016} for the proof of existence of a triangle with given lengths of three
bisectors there was suggested a continuous transformation of a triangle $ABC$ in which 
the bisector's length $l_c$ was approaching to the given length whereas the lengths $l_a$
and $l_b$ were kept constant. Here we shall use the same method but, on the contrary to
Ref.~\cite{os-2016}, such a transformation is considered as a step by step determining of
the elements of the triangles set with prescribed lengths of two bisectors.

For convenience of exposition, this section is subdivided into two parts. In the first part
we consider the set of triangles with $l_a=\mathrm{const}$, $l_b=\mathrm{const}$, and
$0\leq l_c \leq l_a$. In the second part we consider triangles with $l_a=\mathrm{const}$,
$l_b=\mathrm{const}$, and $l_a\leq l_c<\infty$.

\bigskip

\centerline{\large Part 1}

\smallskip

\begin{figure}[ht]
\begin{center}
\includegraphics[width=9cm]{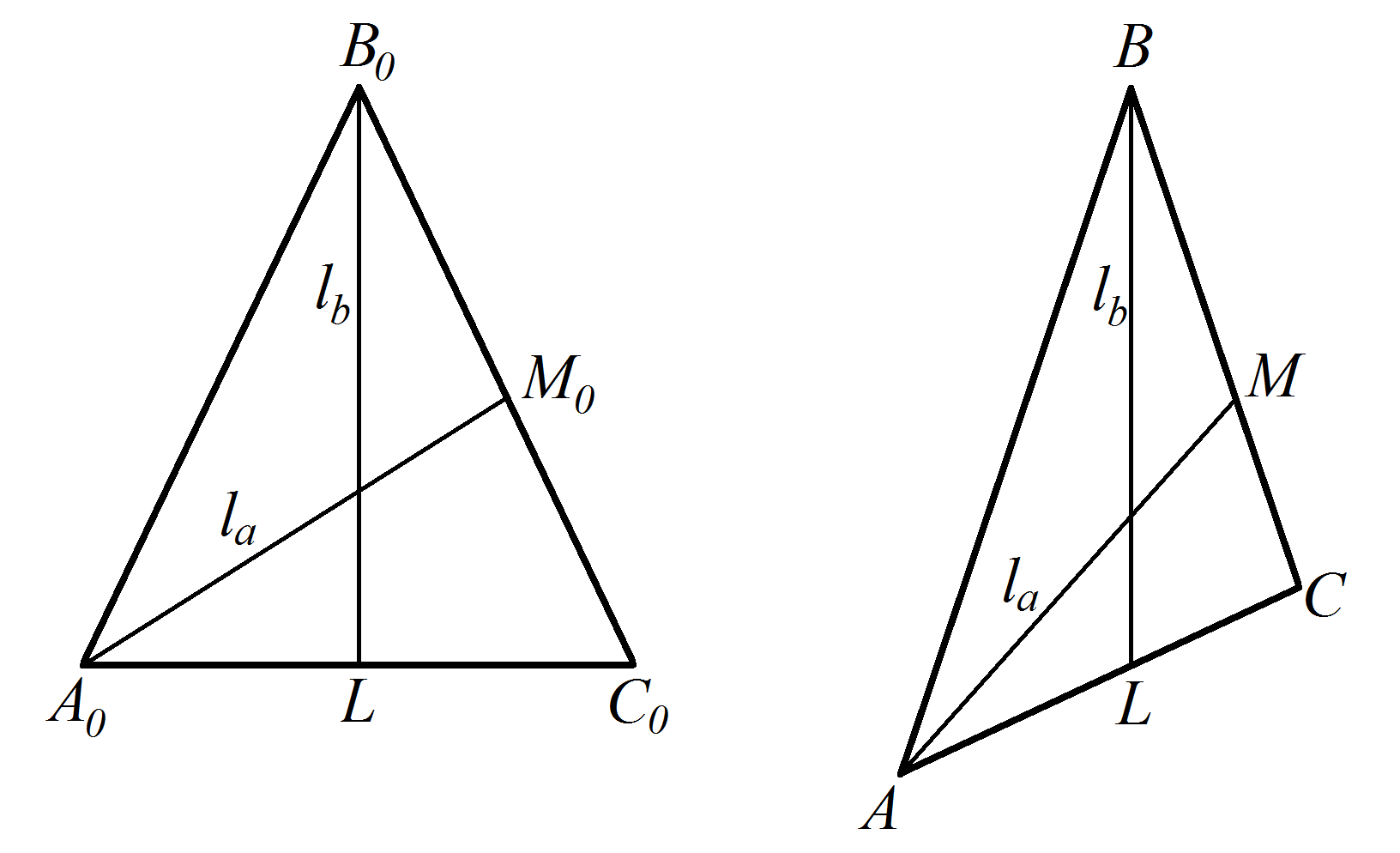}
\put(-200,0){(a)}
\put(-55,0){(b)}
\caption{Construction of triangles $ABC$ with prescribed lengths of the bisectors $l_a=l_2$ and $l_b=l_1$
and the angle $\angle B$ changing within the interval $\eps\leq \angle B\leq\angle B_0$ $(\eps\to0)$.
(a) The triangle $A_0B_0C_0$ is isosceles. (b) Intermediate triangle with $\eps\leq \angle B\leq\angle B_0$.
}
\label{fig4}
\end{center}
\end{figure}

In Fig.~\ref{fig4}(a) the isosceles triangle $ABC$ is depicted with $\angle A=\angle A_0$,
$\angle B=\angle B_0$, $\angle C=\angle C_0$,
$\angle A_0=\angle C_0$, and bisectors $l_a$ and $l_b$ have the prescribed lengths $l_a=l_2$ and $l_b=l_1$.
It is assumed that $l_2\leq l_1$. Notice, that for any triangle it follows from the inequality $l_a\leq l_b$
that $\angle A\geq \angle B$, and vice versa, from $\angle A\geq \angle B$ it follows that $l_a\leq l_b$.
Our task is to get all possible triangles with $l_a=l_2$, $l_b=l_1$, $0\leq l_c\leq l_2$ by means of a continuous
transformation of $ABC$, starting from $A_0B_0C_0$. To this end, we fix the position of the bisector $l_b$
and, consequently, the location of the points $B$ and $L$, but decrease the angle $\angle B$ rotating the
side $AC$ counterclockwise around the point $L$, as is shown in Fig.~\ref{fig4}(b).

It is evident that with increase of the angle $\angle B$ and with fixed value of the angle
$\phi=\angle ALB$ we decrease the sides $AB$ and $AC$, whereas the angle $\angle A$ increases.
From the formula
$$
l_a=\frac{2\cos\frac12\angle A}{1/AB+1/AC}
$$
it follows that $l_a$ increases under this transformation and we can get $l_a<l_2$.
By increasing the rotation angle $\phi$ at  fixed value of $\angle B$ one can increase
$l_a$ up to arbitrary large value (we shall prove this statement in the next section).
Hence, for any value of $\angle B$ there exists only one angle $\phi =\overline{\phi}$
such that $l_a(\overline{\phi})=l_2$. Consequently, by increase of $\angle B$ and simultaneous
increase of $\phi$, we obtain triangles with $l_a=l_2$, $l_b=l_1$. If $\angle B\to0$,
then $\angle A\to0$ and $\angle C\to \pi$ (see Ref.~\cite{os-2016}), and eventually the
triangle $ABC$ degenerates into a segment. The limiting values $l_a$ and $l_b$ will be
equal to $l_a^*=l_2$, $l_b^*=l_1$ and $l_c^*=0$.

Thus, by means of the suggested method we can obtain all possible triangles with $l_a=l_2$,
$l_b=l_1$, and $0\leq l_c\leq l_2$. Let us calculate the limiting values of the sides lengths.
They satisfy the system of equations
$$
l_a^*=\frac{2b^*c^*}{b^*+c^*},\qquad l_b^*=\frac{2a^*c^*}{a^*+c^*},\qquad a^*+b^*=c^*.
$$
This system can be reduced to the quadratic equation
$$
c^{*2} - (l_1 + l_2) c^* +  \frac34 \,l_1\cdot l_2 = 0.
$$
One root only of this equation satisfies the condition $a^{*}<c^{*}$ and $b^*<c^*$.
This root is given by
$$
c^*=\frac12\left(l_1+l_2+\sqrt{l_1^2+l_2^2-l_1l_2}\right),
$$
and we obtain the values of $a^*$ and $b^*$ from the relationships $a^*=c^*/(c^*t_1-1)$,
$b^*=c^*/(c^*t_2-1)$, where $t_1=2/l_1$, $t_2=2/l_2$.

\bigskip

{\it Monotonicity of angles variation}

\smallskip

Since $\angle C=\phi-\angle B/2$, then the angle $\angle C$ increases monotonously with increase of $\phi$
and decrease of $\angle B$. We shall prove now that the angle $\angle A$ decreases monotonously.

\begin{figure}[ht]
\begin{center}
\includegraphics[width=9cm]{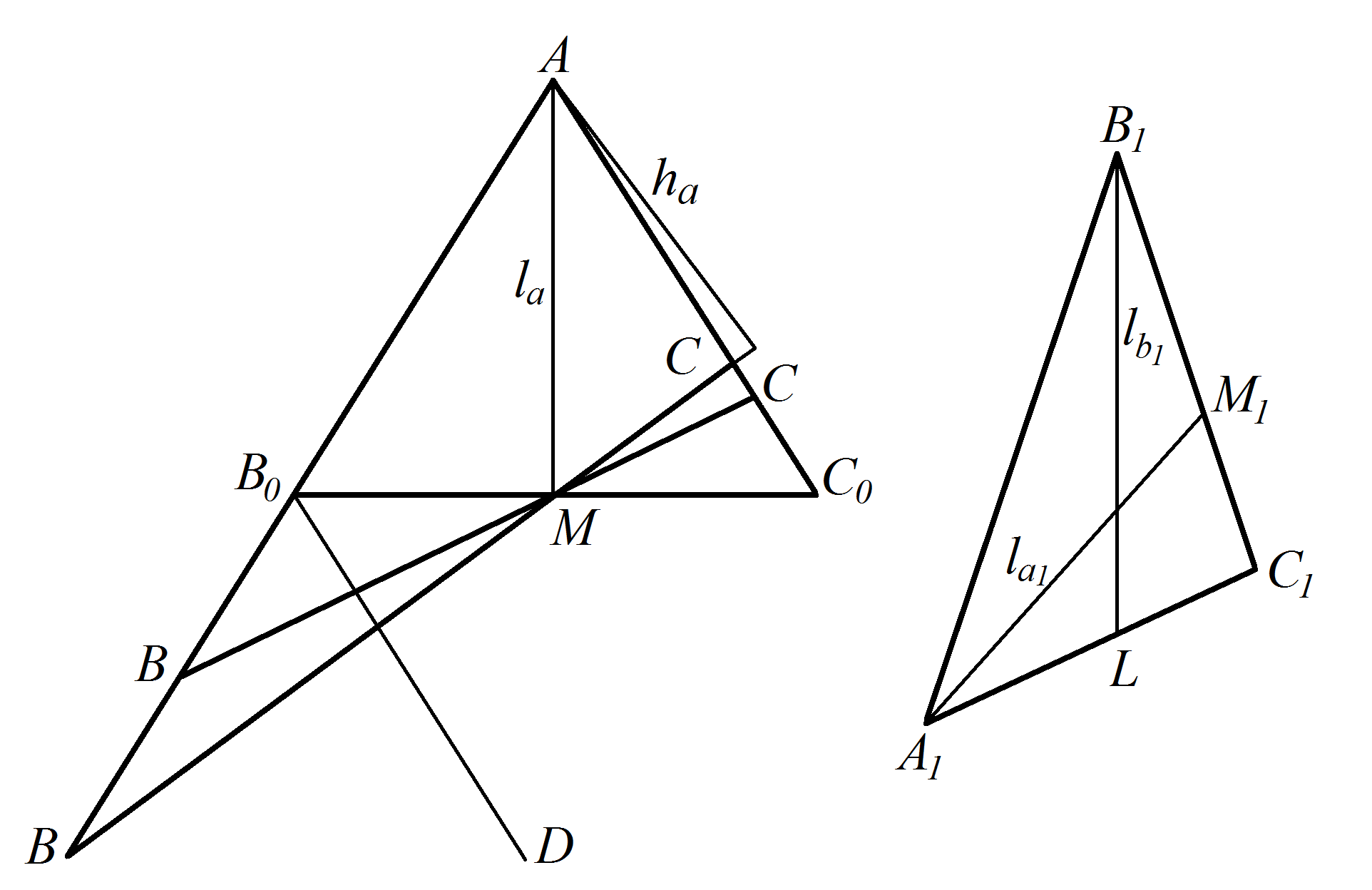}
\put(-210,0){(a)}
\put(-52,0){(b)}
\caption{(a) Sequence of triangles $BAC$ with prescribed values of the angle $\angle A$
and the bisector $l_a=AM$ for the angle $\angle B$ changing from $\angle B=\angle B_0(=\angle C_0)$
to $\angle B=\eps$ $(\eps\to0)$. (b) Triangle $A_1B_1C_1$ which belongs to the set of triangles
$ABC$ ($\angle A=\mathrm{const}$, $l_a=\mathrm{const}$).
}
\label{fig5}
\end{center}
\end{figure}

{\it Lemma 1.} If we draw through the point $M$ taken on the bisector of the angle $\angle A$
the secant line $B_0C_0$ in such a way that $\angle B_0MA =\pi/2$ and after that we rotate this
secant line around the point $M$ counterclockwise, then we get a sequence of triangles $BAC$ (or
a varying triangle $BAC$) such that bisectors of angles of these triangles (or the bisector of the angle
$\angle B$ of the varying triangle) will be the monotonously increasing functions of the rotation
angle (see Fig.~\ref{fig5}(a)).

{\it Proof.}
We draw through the vertex $B_0$ of the triangle $B_0AC_0$ the line $B_0D$ parallel to $AC_0$.
It is clear from Fig.~\ref{fig5}(a) that the area of the triangle $BAC$ increases monotonously
due to an additional part contained inside the angle formed by the lines $B_0B$ and $B_0D$.
But the height $h_a$ drawn from the vertex $A$ down to the side $BC$ decreases monotonously
with growth of the angle $\angle BMA$, because $h_a=AM\sin\angle BMA$ while the angle
$\angle BMA\geq \pi/2$ and it increases during our rotation of the line $BC$. Since the area
of the triangle $BAC$ monotonously increases up to infinity and the height $h_a$ monotonously decreases,
we conclude that the length of $BC$ monotonously increases up to infinity. Obviously, the side $BA$
also monotonously increases and the angle $\angle B$ monotonously decreases down to zero.
Consequently, we infer from the formula
$$
l_b=\frac{2\cos(\angle B/2)}{1/AB+1/BC}
$$
that $l_b$ increases monotonously up to infinity. The proof is finished.

We find from this Lemma that for fixed value of $\angle B$ the bisector $l_b$ increases
monotonously with growth of the rotation angle $\phi$. To prove this statement, it is
enough to exchange the notation of the angles $\angle A$ and $\angle B$.

To prove the monotonous dependence of the angle $\angle A$, we suppose the opposite,
that is that the angle varies non-monotonously. Then there exist at least two triangles
$A_1B_1C_1$ and $A_2B_2C_2$ with $\angle A_1=\angle A_2$, $l_{a1}=l_{a2}=l_2$
$l_{b_1}=l_{b2}=l_1$ and  $\angle B_1\neq\angle B_2$. In Fig.~\ref{fig5}(b) a triangle
$A_1B_1C_1$ is depicted where a bisector $A_1M_1$ of the angle $\angle A_1$ ($A_1M_1=l_{a1}$)
is also shown. Let us estimate the value of the angle $\angle B_1M_1A_1$. Since
$\angle B_1M_1A_1=\tfrac12\angle A_1+\angle C_1$ and $\angle A\geq\angle B$, then
$\angle B_1M_1A_1\geq \tfrac12\angle B_1+\angle C_1=\phi_1\geq\tfrac12\pi$. Similar
estimate takes place for the angle $\angle B_2M_2A_2\geq\tfrac12\pi$. Hence, the triangles
$A_1B_1C_1$ and $A_2B_2C_2$ satisfy the conditions of Lemma~1, and therefore different bisectors
$l_{b1}$ and $l_{b2}$ must correspond to the
different angles $\angle B_1$ and $\angle B_2$,
that is $l_{b1}\neq l_{b2}$, what contradicts to the condition of the
problem that $l_{b1}=l_{b2}=l_1$. Consequently, our supposition about non-monotonicity
of variation of the angle $\angle A$ was wrong.

\bigskip

{\it Monotonicity of sides variation}

\smallskip

Let us prove that $c=AB$ increases monotonously with decrease of the angle $\angle B$.
We denote $\angle A=\al$, $\angle B=\be$, $\angle C=\ga$. Applying sine law to the
triangle $ABC$ we get
$$
c=l_b\frac{\sin(\al+\tfrac12\be)}{\sin\al}=l_a\frac{\sin(\tfrac12\al+\be)}{\sin\be}.
$$
Differentiation of these equalities yields
\begin{equation}\label{eq1}
\begin{split}
  c'(\be)&=l_b\frac{\tfrac12\cos(\al+\tfrac12\be)\sin\al-\al'\sin(\tfrac12\be)}{\sin^2\al}\\
  &=l_a\frac{\tfrac12\al'\cos(\tfrac12\al+\be)\sin\be-\sin(\tfrac12\al)}{\sin^2\be}.
  \end{split}
\end{equation}
According to the imposed conditions we have $l_a\leq l_b$, that is $\al\geq\beta$.
For the isosceles triangle $A_0B_0C_0$ (see Fig.~\ref{fig4}(a)) we get
$\tfrac12\al_0+\beta_0\leq\al_0+\tfrac12\beta_0=\tfrac12\pi$. With decrease of $\beta$
the angle $\al$ decreases also and therefore $\al+\frac12\beta<\tfrac12\pi$ and
$\tfrac12\al+\beta<\tfrac12\pi$. Hence, $\cos(\al+\tfrac12\beta)>0$ and
$\cos(\tfrac12\al+\beta)>0$ (if $\al+\tfrac12\beta=\tfrac12\pi$, then $c'(\beta)<0$).
Let us assume that $c'(\beta)>0$, ($\al+\tfrac12\beta<\tfrac12\pi$). Taking into account
that $\al'(\beta)>0$, we obtain the inequality
$$
\frac{\frac12\cos(\al+\tfrac12\beta)\sin\al}{\sin\tfrac12\beta}>\al'>
\frac{\sin\tfrac12\al}{\frac12\cos(\tfrac12\al+\beta)\sin\beta}.
$$
This yields
$$
\frac{\frac12\cos(\al+\tfrac12\beta)\sin\al}{\sin\tfrac12\beta}>
\frac{\sin\tfrac12\al}{\frac12\cos(\tfrac12\al+\beta)\sin\beta}
$$
and after simplifications we get
$$
\cos(\al+\tfrac12\beta)\cos(\tfrac12\al+\beta)\cos\tfrac12\al\cos\tfrac12\beta>1,
$$
what is impossible. Hence, $c'(\beta)<0$ and $c(\beta)$ increases monotonously with growth
of $\beta$. In addition, we have proved the inequality
\begin{equation}\label{eq2}
  \frac{\frac12\cos(\al+\tfrac12\beta)\sin\al}{\sin\tfrac12\beta}<\al'<
\frac{\sin\tfrac12\al}{\frac12\cos(\tfrac12\al+\beta)\sin\beta}.
\end{equation}

Now let us prove that the sides $a$ and $b$ are monotonously decreasing.
Indeed, from the relations
\begin{equation}\label{eq3}
  l_a=\frac{2\cos(\tfrac12\al)}{1/b+1/c},\qquad l_b=\frac{2\cos(\tfrac12\beta)}{1/a+1/c}
\end{equation}
it follows that for keeping the constant values of $l_a$ and $l_b$ it is necessary to have
$a$ and $b$ monotonously decreasing.

\bigskip

{\it Monotonicity of heights variation}

\smallskip

Let us prove that all heights decrease monotonously down to zero with decrease of
$\beta$ down to $\beta=0$. We have
\begin{equation}\label{eq4}
  \begin{split}
  & h_a=c\sin\beta=l_a\sin(\tfrac12\al+\beta),\\
  & h_b=c\sin\al=l_a\sin(\al+\tfrac12\beta),\\
  & h_c=a\sin\beta.
  \end{split}
\end{equation}
Since $\al$ and $\beta$ decrease monotonously down to zero, the heights $h_a$ and $h_b$
also decrease monotonously down to zero. But for $\al+\tfrac12\beta=\tfrac12\pi$ the
height $h_b$ has the maximal value $h_b=l_b$. The height $h_c$ decreases also down to
zero because $\sin\beta$ decreases down to zero and the side $a$ decreases monotonously
down to $a=a^*$. If $l_a=l_b$, then the triangle $A_0B_0C_0$ is equilateral
($\angle A_0=\angle B_0=\angle C_0=\pi/3$). In this case $\tfrac12\al_0+\beta_0=\tfrac12\pi$
and the height $h_a$ reaches its maximal value $h_a=l_a$.

\bigskip

{\it Monotonicity of bisector $l_c$ variation}

\smallskip

We have shown above that with decrease of the angle $\angle B$ from $\angle B=\angle B_0$
down to $\angle B=0$, the sides $a$ and $b$ decrease monotonously down to $a=a^*$
and $b=b^*$, whereas $\angle C$ increases monotonously up to $\angle C=\pi$. Therefore the
bisector $l_c$ decreases monotonously from the value $l_c=l_a$ down to the value $l_c=0$.

\bigskip

{\it Monotonicity of medians variation}

\smallskip

(i) Let us show that the median $m_c$ decreases monotonously with decrease of the angle $\angle B$.
To this end, we shall use the known relationship between the median and the triangle's sides,
$$
4m_c^2=2a^2+2b^2-c^2.
$$
Differentiation with respect to $\beta$ gives
$$
4m_cm_c'=2aa'+2bb'-cc'.
$$
Since $a'>0$, $b'>0$, and $c'<0$, then $m_c'>0$. Hence, when $\beta$ decreases,
the median $m_c$ also decreases and at $\beta=0$ it takes the value
$m_c=m_c^*=\frac12(a^*-b^*)$.

(ii) Now we shall find how the median $m_b$ changes with decrease of the angle $\angle B$.
We take the expression for $m_b$ in terms of the triangle's sides and transform it with the use
of the law of cosines,
\begin{equation}\label{eq5}
  \begin{split}
  4m_b^2&=2c^2+2a^2-b^2=c^2+a^2+2ac\cos\beta=\\
  &=(c-a)^2+2ac(1+\cos\beta)=(c-a)^2+4ac\cos^2\left(\tfrac12\beta\right).
  \end{split}
\end{equation}
We take the derivative of this expression to get
$$
4m_bm_b'=(c-a)(c'-a')+2ac\cos^2\left(\tfrac12\beta\right)
\left[\ln\left(ac \cos^2\left(\tfrac12\beta\right)\right)\right]'.
$$
We denote
$$
u_1=(c-a)(c'-a'),\qquad v_1=2ac\cos^2\left(\tfrac12\beta\right)
\left[\ln\left(ac \cos^2\left(\tfrac12\beta\right)\right)\right]'.
$$
Since $c'<0$, $a'>0$, then $c'-a'<0$, and, hence, $u_1<0$, if $c>a$, $u_1>0$ if $c<a$, and
$u_1=0$, if $c=a$. Now we transform $v_1$,
\begin{equation}\label{eq6}
  \begin{split}
  v_1&=2ac\cos^2\left(\tfrac12\beta\right)\left[\ln\left(\frac{h_b}{\sin(\alpha+\beta)}
  \cdot\frac{h_b}{\sin\alpha}\cdot \cos^2\left(\tfrac12\beta\right)\right)\right]=\\
  &=2ac\cos^2\left(\tfrac12\beta\right)
  \left[2\frac{h_b'}{h_b}-\frac{\cos(\alpha+\beta)}{\sin(\al+\beta)}(\alpha'+1)
  -\frac{\cos\alpha}{\sin\al}\,\alpha'-
  \frac{\sin\tfrac12\beta}{\cos\tfrac12\beta}\right].
  \end{split}
\end{equation}
We substitute here $-h_b'/h_b=\cot(\alpha+\beta/2)(\alpha'+1/2)$ and obtain after some
simplifications the expression
$$
v_1=\frac{ac\sin\beta\left[-\alpha'\sin\beta+\cos(\alpha+\beta)\sin\alpha\right]}
{\sin(\alpha+\beta)\sin\alpha}\,\cot(\alpha+\beta/2).
$$
Let us show that
\begin{equation}\label{eq7}
  -\alpha'\sin\beta+\cos(\alpha+\beta)\sin\alpha<0.
\end{equation}
Indeed, for $\alpha+\beta\geq\tfrac12\pi$ the validity of (\ref{eq7}) is obvious.
For $\alpha+\beta<\tfrac12\pi$ we get from (\ref{eq7})
${\alpha'>\cos(\alpha+\beta)\sin\alpha}/{\sin\beta}$. It was shown above (see (\ref{eq2}))
that if $\alpha+\beta<\tfrac12\pi$, then
$$
\alpha'>\frac{\cos(\alpha+\beta/2)\sin\alpha}{2\sin(\beta/2)}.
$$
We obviously have
$$
\frac{\cos(\alpha+\beta/2)\sin\alpha}{2\sin(\beta/2)}>\frac{\cos(\alpha+\beta)\sin\alpha}{sin\beta},
$$
and this inequality can be transformed to
$$
\sin(\alpha+\beta/2)\sin(\beta/2)>0.
$$
This proves the validity of (\ref{eq7}). The other factors in expression for $v_1$ are positive
or vanishing, therefore $v_1\leq0$, and since $a\leq c$, then we have $u_1\leq0$. If $a=c$ and $\alpha+\beta/2=\pi/2$,
then $u_1=v_1=0$. In this case $m_b$ takes its minimal value $m_b=l_b$. With decrease of $\beta$
down to $\beta=0$ the median $m_b$ increases monotonously up to the value
$m_b=m_b^*=a^*+\tfrac12b^*$.

(iii) Now let us consider variation of the median $m_a$. We have
\begin{equation}\label{eq8}
4m_a^2=2b^2+2c^2-a^2=(c-b)^2+4bc\cos^2(\alpha/2).
\end{equation}
Differentiation yields
$$
4m_am_a'=(c-b)(c'-b')+2bc\cos^2(\alpha/2)\left[\ln\left(bc\cos^2(\alpha/2)\right)\right]'.
$$
We denote
$$
u_2=(c-b)(c'-b'),\qquad v_2=2bc\cos^2(\alpha/2)\left[\ln\left(bc\cos^2(\alpha/2)\right)\right]'.
$$
Since $b'>0$, $c'<0$, then $c'-b'<0$ and therefore $u_2<0$ if $c>b$, $u_2>0$ if $c<b$, or $u_2=0$ if $c=b$.
Expressing $b$ and $c$ as functions of $\alpha$ and $\beta$, we obtain
\begin{eqnarray}
  b &=& \frac{h_a}{\sin(\alpha+\beta)}=\frac{l_a\sin((\alpha+\beta)/2)}{\sin(\alpha+\beta)}, \\
  c &=& \frac{h_a}{\sin\beta}=\frac{l_a\sin((\alpha+\beta)/2)}{\sin\beta}.
\end{eqnarray}
Substitution of these expressions into formula for $v_2$ and comparison with $v_1(\alpha,\beta)$
show that
$$
\left[\ln\left(ac\cos^2(\beta/2)\right)\right]_{\beta}'=q(\alpha,\beta),\quad\text{and}\quad
\left[\ln\left(bc\cos^2(\alpha/2)\right)\right]_{\alpha}'=q(\beta,\alpha).
$$
Consequently,
\begin{equation}\nonumber
  \begin{split}
  v_2&=2bc\cos^2(\alpha/2)\left[\ln\left(bc\cos^2(\alpha/2)\right)\right]_{\alpha}'\alpha_{\beta}'\\
  &=\alpha'\frac{bc\sin\alpha\left[-\beta'\sin\alpha+\cos(\alpha+\beta)\sin\beta\right]}
  {\sin(\alpha+\beta)\sin\beta}\,\cot(\alpha/2+\beta).
  \end{split}
\end{equation}
Let us consider the product
$$
p=\alpha'\left[-\beta'\sin\alpha+\cos(\alpha+\beta)\sin\beta\right]=
-\sin\alpha+\alpha'\cos(\alpha+\beta)\sin\beta
$$
and show that
\begin{equation}\label{eq11}
  p<0.
\end{equation}
For $\alpha+\beta\geq\pi/2$ the validity of this inequality is obvious. Let  $\alpha+\beta<\pi/2$ and then
we have to show that $\alpha'<\sin\alpha/(\cos(\alpha+\beta)\sin\beta)$. It was shown above (see (\ref{eq2}))
that if $\alpha+\beta/2<\pi/2$, then $\alpha'<2\sin(\alpha/2)/(\cos(\alpha/2+\beta)\sin\beta)$. Since
$$
\frac{2\sin\alpha}{\cos(\alpha/2+\beta)}<\frac{\sin\alpha}{\cos(\alpha+\beta)\sin\beta},
$$
this inequality can be transformed to
$$
\sin(\alpha/2+\beta)\sin(\alpha/2)>0.
$$
The other factors are positive (except for the case $l_a=l_b$), hence $v_2<0$, and since $c>b$,
we get also $u_2<0$. Therefore the median $m_a$ increases monotonously with decrease of $\beta$
up to the value $m_a=m_a^*=b^*+a^*/2$. If the triangle $A_0B_0C_0$ is equilateral, then
$b_0=c_0$, $\alpha_0/2+\beta_0=\pi/2$, $u_2(b_0,c_0)=0$, $u_2(\alpha_0,\beta_0)=0$ and the
median takes the minimal value $m_a=l_a$.

\bigskip
\bigskip
\bigskip

\centerline{\large Part 2}

\smallskip

\begin{figure}[ht]
\begin{center}
\includegraphics[width=9cm]{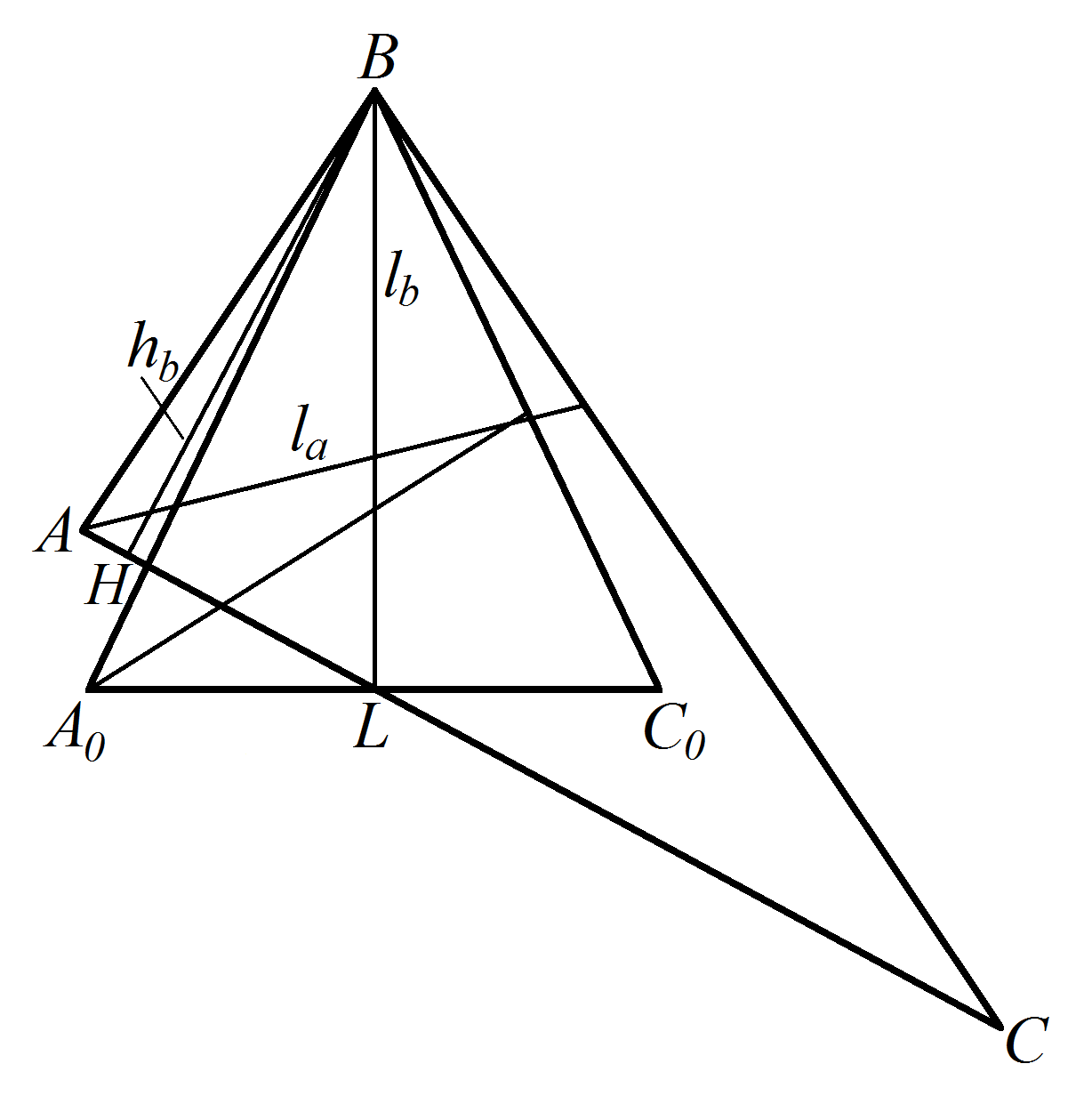}
\caption{Formation of triangles with fixed values of $l_a=l_2$ and $l_b=l_1$ by variation of
the angle $\angle B$ starting from the value $\angle B=\angle B_0$ $(\angle A_0=\angle C_0)$.
}
\label{fig6}
\end{center}
\end{figure}

Let us consider the set of triangles $ABC$ with $l_a=\mathrm{const}=l_2$, $l_b=\mathrm{const}=l_1$,
and $l_c\in \{l_2\leq l_c<\infty\}$. These triangles can be obtained by the method of continuous
transformation which was used in Part~1 after some its modification. As earlier, this transformation starts
from an isosceles triangle $A_0B_0C_0$ shown in Fig.~\ref{fig4}(a)), but now we increase $\angle B$ and
rotate the side $AC$ clockwise around the point $L$ in such a way, that the location of the bisector $l_b=BL$
does not change in the plane of the triangle (see Fig.~\ref{fig6}).

\begin{figure}[ht]
\begin{center}
\includegraphics[width=9cm]{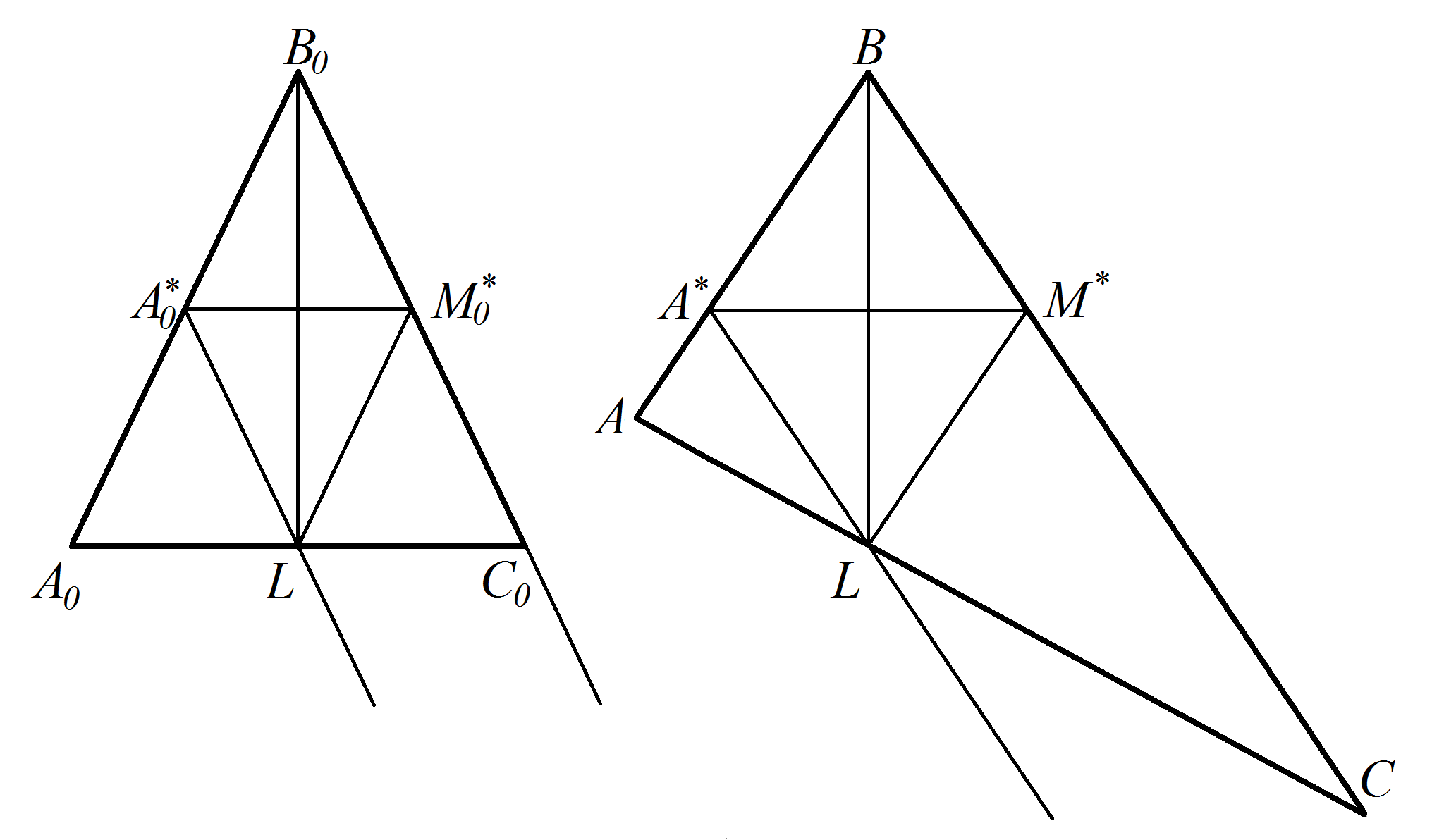}
\put(-210,0){(a)}
\put(-105,0){(b)}
\caption{(a) Limiting position of the bisector $l_a=A_0^*M_0^*$ at the maximal rotation of the side
$A_0C_0$ and fixed value of the angle $\angle B=\angle B_0$.
(b) Limiting position of $l_a=A^*M^*$ at $\angle B_0<\angle B\leq\angle B_{\text{max}}$.
}
\label{fig7}
\end{center}
\end{figure}

Let us show that the increase of the angle $\angle B$ leads to the growth of $l_a$, so that
the equality $l_a=l_2$ stops to hold, whereas rotation of $AC$ around the point $L$, on the
contrary, decreases its length down to $l_a<l_2$. As a result of both transformation the
length of $l_a$ can be restored and kept constant while the angle $\angle B$ increases.
We denote the rotation angle as $\angle BLC=\phi_1$. It is evident, that as long as the
height $h_b$ is located inside the angle $\angle B$, the increase of $\angle B$ leads to
the increase of the sides $AB$ and $AC$ and to the decrease of $\angle A$, so $l_a$ increases.
However, with further increase of $\phi_1$ the height $h_b$ exits the interior of the angle $\angle B$ 
and the point $H$ becomes located on the side $AC$ prolonged outside $\angle B$ (see Fig.~\ref{fig6}).
In this case the side $AB$ decreases with increase of $\angle B$ and the special study is
necessary for proving of growth of $l_a$. We denote $\angle HAB=\psi$. If we increase
$\angle B$ without changing $\phi_1$, then $h_b$ does not change either. Substitution of
$c=h_b/\sin\psi$ into the formula for $l_a$ yields
\begin{equation}\label{eq12}
  l_a=\frac{2h_bb/\sin\psi}{b+c}\cos\left[\tfrac12(\pi-\psi)\right]=\frac{{h_b}/{\cos(\psi/2)}}{c/b+1}.
\end{equation}
Since with the increase of the angle $\angle B$ both the angle $\psi$ and the side $b$ are
increasing, the side $c$ is decreasing, then it follows from Eq.~(\ref{eq12}) that the length $l_a$ 
increases and eventually we get $l_a>l_2$. Let us consider, how the increase of $\phi_1$ influences 
on the length $l_a$. We take an isosceles triangle $A_0B_0C_0$ and rotate $A_0C_0$ around the point 
$L$ by the angle $\phi_1=\phi_{10}^*$ such that the straight line $A_0^*L$ becomes parallel to the
line $B_0C_0$. Now we draw the bisector $A_0^*M_0^*=l_{a0}^{**}$ of the angle $\angle A_0^*$
(see Fig.~\ref{fig7}(a)).

The quadrilateral $A_0^*B_0M_0^*L$ is rhombus. Its diagonal $A_0^*M_0^*$ is equal to
$A_0^*M_0^*=l_{a0}^{**}=\tfrac12A_0C_0$. Let us show that $l_{a0}^{**}<l_2$.

{\it Lemma 2.} In any isosceles triangle $ABC$ with base $AC=b$ and angle $\angle A=\al$
the inequality
$$
\tfrac23b\leq l_a\leq \sqrt{2}\,b
$$
holds; besides that, $l_a=\tfrac23b$
if $\al=0$ and $l_a=\sqrt{2}\,b$ if $\al=\tfrac12\pi$.

{\it Proof.} We take an arbitrary isosceles triangle $ABC$ ($\angle A=\angle C$) with the
base $AC=b=\mathrm{const}$ and find how the ratio $l_a/b$ changes as a function of $\al$.
The derivative of $l_a/b$ is equal to
$$
\left(\frac{l_a}b\right)'=\left(\frac{\sin\al}{\sin\frac{3\al}2)}\right)'=
\frac{\sin\frac{\al}2\left(1-\cos\frac{\al}2\cos\frac{3\al}2\right)}{\sin^2\frac{3\al}2}
>0\quad\text{or}\quad =0\quad\text{if}\quad \al=0.
$$
Consequently,
$$
\mathrm{min}\,\frac{l_a}{b}=\lim_{\al\to0}\frac{\sin\al}{\sin(3\al/2)}=\frac23\qquad\text{and}
\qquad
\mathrm{max}\,\frac{l_a}{b}=\lim_{\al\to\tfrac{\pi}2}\frac{\sin\al}{\sin(3\al/2)}=\sqrt{2}.
$$
This finishes the proof.

Since in the triangle $A_0B_0C_0$ we have $l_{a0}=l_2$ and $\tfrac12A_0C_0=\tfrac12b_0=l_{a0}^{**}$,
then Lemma 2 yields
$$
\frac{\sqrt{2}}4l_2<l_{a0}^{**}<\frac34l_2,
$$
and because $\mathrm{min}\,\al_0=\tfrac13\pi$ (rather than zero), we get
$$
\frac{\sqrt{2}}4l_2<l_{a0}^{**}<\frac{\sqrt{3}}3l_2.
$$

The importance of the result obtained here consists of the possibility to make a statement that
one can increase the angle $\angle B$ from the value $\angle B=\angle B_0$ up to some
maximal value $\angle B=\angle B_{\text{max}}$ with fulfilling the condition $A^*M^*\leq l_2$,
where $A^*M^*$ is the diagonal of the rhombus $A^*BM^*L$ (see Fig.~\ref{fig7}(b)).

To justify this statement, we fix the value of $\angle B$ at some intermediate level between
$\angle B_0$ and $\angle B_{\text{max}}$ and rotate $AC$ around the point $L$ up to the
maximal angle $\phi_1^*$. Now we draw the bisector of the angle $\angle A^*$ and denote
as $M^*$ its intersection point with $BC$ (see Fig.~\ref{fig7}(b)). The quadrilateral
$A^*BM^*L$ is the rhombus which differs from the rhombus $A_0^*B_0M_0^*L$ by the angle
$\angle B>\angle B_0$ and by the length of the diagonal $A^*M^*=l_a^{**}$ (the length
$l_a^{**}>l_{a0}^{**}$ and it increases with growth of $\angle B$). During the rotation of
$AC$ around the point $L$ the bisector $l_a$ decreases monotonously from $l_a>l_2$ down to
the value $l_a=l_a^{**}$ (see Lemma 3 below).
Consequently, as long as $l_a^{**}<l_2$, for each $\angle B\geq \angle B_0$ the exists a
unique rotation angle $\overline{\phi}_1<\phi_1^*$, which depends on $\angle B$, such that
$l_a(\overline{\phi}_1)=l_2$. Hence, as long as $l_a^{**}<l_2$, for each
$\angle B\geq \angle B_0$ there exists a unique triangle $ABC$ in which $l_a=l_b$ and $l_b=l_1$.
As soon as we have $\angle B=\angle B_{\text{max}}=2\arctan(l_2/l_1)$, $l_a^{**}=l_2$ and
$\overline{\phi}_1=\phi_1^*$, the triangle $ABC$ stops to exist converting into a biangle
with parallel sides $A^*L$ and $BC$ (see Fig.~\ref{fig7}(b)). Therefore for
$\angle B\to\angle B_{\text{max}}$ we have $AC\to\infty$, $BC\to\infty$, $\angle C\to0$,
and $l_c\to\infty$. Thus, proceeding in this way with varying $\angle B$ within the interval
$\angle B_0\leq \angle B\leq\angle B_{\text{max}}$, we obtain all triangles with
$l_2\leq l_c<\infty$ and, with account of the result obtained in Part~1, we obtain the whole
set of triangles $ABC$ with $l_a=l_2$, $l_b=l_1$ and $0\leq l_c<\infty$.

One should note that $\angle B_{\text{max}}\leq\tfrac12\pi$ (the equality corresponds to the
equilateral triangle $A_0B_0C_0$). The limiting value of the angle
$\angle A=\angle A_{\text{max}}=\pi-\angle B_{\text{max}}\geq\frac12\pi$

\bigskip

{\it Monotonicity of angles variation}

\smallskip

It is evident that with monotonous increase of the angles $\angle B$ and $\phi_1=\angle BLC$,
the angle $\angle C$ decreases monotonously down to the value $\angle C=0$.

\begin{figure}[ht]
\begin{center}
\includegraphics[width=9cm]{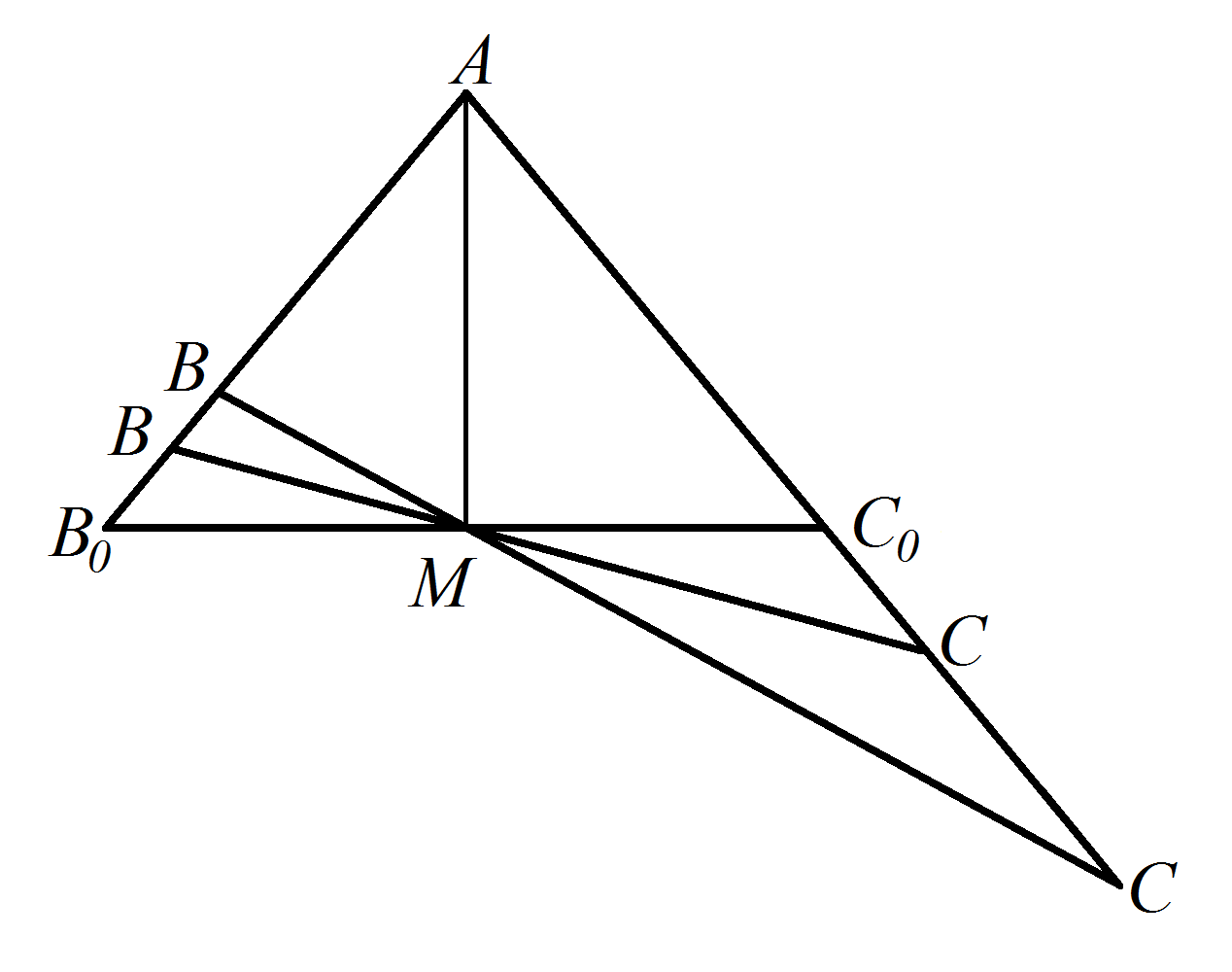}
\caption{The sequence of triangles $BAC$ with fixed $\angle A$ and $l_a=AM$. The angle
$\angle B$ increases monotonously from $\angle B=\angle B_0\,(=\angle C_0)$ up to
$\angle B=\pi-\angle A$.
}
\label{fig8}
\end{center}
\end{figure}

Let us prove that $\angle A$ increases monotonously in close analogy with the proof presented in
Part 1. We consider a monotonous sequence of triangles $BAC$ which have $\angle A=\mathrm{const}$,
$l_a=AM=\mathrm{const}$ and $\angle AMC$ increases from $\angle AMC=\tfrac12\pi$ up to
$\angle A=\pi-\tfrac12\angle A$ (see Fig.~\ref{fig8}).

Now we show that the bisector of the angle $\angle B$ decreases monotonously with growth of
$\angle AMC$. We denote $\angle A=\al$, $\angle B=\beta$, $l_a=l$, and we have
$$
l_b=\frac{2ac}{a+c}\cos\left(\frac{\beta}2\right),\quad c=l\,\frac{\sin(\al/2+\beta)}{\sin\beta},
\quad
a=c\,\frac{\sin\al}{\sin(\al+\beta)}.
$$
From these formulas we get
$$
l_b=\frac{2c}{a/c+1}\cos\left(\frac{\beta}2\right)=
\frac{2l\sin(\al/2+\beta)\cos(\beta/2)}{\sin\beta(\frac{\sin(\al+\beta)}{\sin\al}+1)}=
\frac{l\sin\al\sin(\al/2+\beta)}{\sin\frac{\beta}2(\sin(\al+\beta)+\sin\al)}.
$$
The derivative $l_b'(\beta)$ for constant $\al$ and $l$ is equal to
\begin{equation}\nonumber
  \begin{split}
  \frac{l_b'(\beta)}{l\sin\al}=&\Big\{\cos(\frac{\al}2+\beta)\sin\frac{\beta}2[\sin(\al+\beta)+\sin\al]\\
  &-\sin(\frac{\al}2+\beta)
  \{\tfrac12\cos\frac{\beta}2[\sin(\al+\beta)+\sin\al]+\sin\frac{\beta}2\cos(\al+\beta)\}\Big\}\\
  &\times\left[\sin^2\frac{\beta}2[\sin(\al+\beta)+\sin\al]^2\right]^{-1}.
  \end{split}
\end{equation}
We find
\begin{equation}\nonumber
  \begin{split}
  & f_1=-(1-\cos\al)\cos(\frac{\al+\beta}2)\leq0\qquad (=0\quad\text{if}\quad \al+\beta=\pi),\\
  & f_2=-\sin(\al/2+\beta)\cos\frac{\beta}2\sin(\al+\beta)\leq0\qquad (=0\quad\text{if}\quad \al+\beta=\pi),\\
  & f_3=\cos(\al/2+\beta)\sin\frac{\beta}2\sin\al\leq0\qquad (\text{since}\quad \al/2+\beta\geq\pi/2).
  \end{split}
\end{equation}
We notice that $f_1,f_2,f_3$ do not vanish simultaneously. It is not difficult to check the validity of
the equality
\begin{equation}\label{eq13}
  l_b'(\beta)=\frac{l}2\sin\al\frac{f_1+f_2+f_3}{\sin^2\frac{\beta}2[\sin(\al+\beta)+\sin\al]^2}.
\end{equation}
It follows from (\ref{eq13}) that $l_b'(\beta)<0$ and, consequently, the bisector $l_b$ decreases
monotonously with growth of $\beta$.

Interchanging the notation of the angles $\angle A$ and $\angle B$, we formulate the obtained result
as the following

{\it Lemma 3.} With a clockwise increase of the rotation angle of the side $AC$ of the triangle $ABC$
around the point $L$ at the fixed value of the angle $\angle B$, the bisector of the angle
$\angle A$ decreases down to the some value which depends on $\angle B$.

Turning to the proof of the main statement, we suppose that $\angle A$ changes non-monotonously.
Then there exist at least two triangles $A_1B_1C_1$ and $A_2B_2C_2$ with equal angles
$\angle A_1$ and $\angle A_2$, equal bisectors of these angles and different angles $\angle B_1$
and $\angle B_2$. We superimpose the angles $\angle A_1$ and $\angle A_2$ and their bisectors
of these triangles. According to the obtained above result, the bisectors of the angles
$\angle B_1$ and $\angle B_2$ must have different lengths what contradicts to the condition
$l_b=\textrm{const}$. Thus, our supposition is wrong and the angle $\angle A$ changes
(increases) monotonously.

\bigskip

{\it Monotonicity of sides variation}

\smallskip

From (\ref{eq1}) we see that if $\al+\be/2\geq \pi/2$ ($\cos(\al+\be/2)\leq0$), then $c'(\be)<0$.
Consequently, with growth of $\be$ from $\be=\be^{**}=2\arctan(l_2/l_1)$ the side $c$ decreases
monotonously down to the value $c=c^{**}=\tfrac12\sqrt{l_1^2+l_2^2}$, and from (\ref{eq3}) it
follows that $a$ and $b$ increase monotonously up to infinity.

\bigskip

{\it Specific features of heights variation}

\smallskip

Let us turn to the formulas (\ref{eq4}).

(i) $h_a=l_a\sin(\al/2+\be)$.

If $l_1=l_2$ (the triangle $A_0B_0C_0$ is isosceles), then $\al_0/2+\be_0=\pi/2$ and $h_a$ takes its
maximal value $h_a=h_a(\al_0,\be_0)=l_a$. With increase of $\al$ and $\be$ up to their maximal
values $\al=\al^{**}=\pi-\be^{**}$ and $\be=\be^{**}$ the height $h_a$ decreases monotonously
down the value $h_a=h_a^{**}=l_1l_2/\sqrt{l_1^2+l_2^2}$. If $l_2<l_1$, then
$\al_0/2+\be_0<\al_0+\be_0/2=\pi/2$, but at $\al+\be=\pi$ we have
$\al/2+\be=(\al+\be)/2)+\be/2=\pi/2+\be/2>\pi/2$. Consequently, there exist such $\al=\overline{\al}$
and $\be=\overline{\be}$, that $\overline{\al}/2+\overline{\be}=\pi/2$. Therefore
$h_a=h_a(\overline{\al},\overline{\be})=l_a$ is the maximal value of $h_a$.

(ii) $h_b=l_b\sin(\al+\be/2)$

Since in the triangle $A_0B_0C_0$ the equality $\al_0+\be/2=\pi/2$ is fulfilled, then
$h_b=h_b(\al_0,\be_0)=l_b$ is the maximal value of $h_b$. With growth of $\al$ and $\be$
up to their maximal values $\al=\al^{**}$ and $\be=\be^{**}$, the value of $h_b$
decreases monotonously down its minimum $h_b=h_b^{**}=l_1l_2/\sqrt{l_1^2+l_2^2}$.

(iii) $h_c=a\sin\be$

Since both $a$ and $\sin\be$ increase monotonously with growth of $\be$
($\be\to\be^{**}\leq\pi/2$ and $a\to\infty$), then $h_c$ increases monotonously up to
infinity.

\bigskip

{\it Monotonous growth of the bisector $l_c$}

\smallskip

As was shown above, with growth of $\be$ from $\be=\be_0$ up to $\be=\be^{**}$,
the sides $a$ and $b$ increase monotonously up to infinity, whereas the angle
$\angle C$ decreases monotonously down to zero. Consequently, the bisector $l_c$
increases monotonously from $l_c=l_2$ up to $l_c=\infty$.

\bigskip

{\it Specific features of medians variation}

\smallskip

To analyze the specific features of variation of medians for $\be$ increasing in the interval
under consideration, one has to use the transformations introduced in the Part~1. Taking into
account that $\al+\be/2\geq\pi/2$ and $\al/2+\be\leq\al+\be/2$ (the equality sign corresponds to
$l_1=l_2$), we arrive at the following conclusions.

(i) The median $m_c$ increases monotonously up to $m_c=\infty$, because $l_c\to\infty$ and $m_c\geq l_c$.

(ii) The median $m_b$ increases monotonously from its minimal value $m_b=l_b=l_1$ up to
$m_b=\infty$. Indeed, from the formula $4m_b^2=(a-c)^2+4ac\cos^2(\be/2)$ with $a\to\infty$ and
$c\to c^{**}$ we obtain at once that $m_b\to\infty$ for $\be\to\be^{**}$.

(iii) If $l_2<l_1$, then $\al/2+\be<\al+\be/2$ and, consequently, when $\al+\be/2=\pi/2$ we get
$\al/2+\be<\pi/2$. But with increase of $\al$ and $\be$ the sum $\al/2+\be$ increases, too,
and at $\al=\overline{\al}$, $\be=\overline{\be}$ (see above) we obtain $\overline{al}/2+\overline{\be}=\pi/2$,
so the median $m_a$ decreases monotonously down to its minimal value
$m_a=m_a(\overline{\al},\overline{\be})=l_a=l_2$. With further increase of $\al$ and $\be$ the median
$m_b$ increases monotonously up to $m_a=\infty$, because in the relationship
$4m_a^2=(b-c)^2+4bc\cos^2(\al/2)$ we have $b\to\infty$ and $c\to c^{**}$. If $l_1=l_2$, then
$m_a$ changes in the same way as $m_b$.

\section{Results and conclusions}

The results obtained here permit us to make the certain statements about existence of triangles with
prescribed two bisectors and one third element which can be taken as one of the angles, the sides,
the heights or the medians, or the third bisector. Let $l_a=l_2$ and $l_b=l_1$, where $l_1$ and $l_2$ are
given lengths satisfying the condition $l_2\leq l_1$.

\bigskip

\begin{enumerate}
  \item Let the angle $\angle A$ be given. If its value belongs to the interval
  $0<\angle A<\angle A_{\text{max}}$, then there exists a unique triangle with prescribed values of
  $l_a$, $l_b$ and $\angle A$.
  \item Let the angle $\angle B$ be given. If its value belongs to the interval
  $0<\angle B<\angle B_{\text{max}}$, then there exists a unique triangle with prescribed values of
  $l_a$, $l_b$ and $\angle B$.
  \item Let the angle $\angle C$ be given. Then for any $\angle C$ from the interval
  $0<\angle C<\pi$ there exists a unique triangle with prescribed values of
  $l_a$, $l_b$ and $\angle C$.
  \item Let the side $a$ be given. Then for any $a>a^*$ there exists a unique triangle with
  prescribed $l_a$, $l_b$ and $a$.
  \item Let the side $b$ be given. Then for any $b>b^*$ there exists a unique triangle with
  prescribed $l_a$, $l_b$ and $b$.
  \item Let the side $c$ be given. If its length belongs to the interval $c^{**}< c< c^*$,
  then there exists a unique triangle with prescribed $l_a$, $l_b$ and $c$.
  \item Let the height $h_a$ be given. If its length belongs to the interval $h_a^{**}< h_a< l_b$,
  then there exist two triangles with prescribed $l_a$, $l_b$ and $h_a$. If $h_a$ belongs to
  the interval $0<h_a\leq h_a^{**}$ or $h_a=l_a$, then there exists a unique triangle
  with prescribed $l_a$, $l_b$ and $h_a$.
  \item Let the height $h_b$ be given. If its length belongs to the interval $h_b^{**}< h_a< l_a$,
  then there exist two triangles with prescribed $l_a$, $l_b$ and $h_b$. If $h_b$ belongs to
  the interval $0<h_b\leq h_b^{**}$ or $h_b=l_b$, then there exists a unique triangle
  with prescribed $l_a$, $l_b$ and $h_b$.
  \item Let the height $h_c$ be given. Since $h_c$ varies continuously from zero to infinity, then
  for any $h_c$ from the interval $0<h_c<\infty$ there exists a unique triangle
  with prescribed $l_a$, $l_b$ and $h_c$.
  \item Let the median $m_a$ be given. If its length belongs to the interval $l_a<m_a<m_a^{*}$,
  then there exist two triangles with prescribed $l_a$, $l_b$ and $m_a$. If $m_a$ belongs to
  the interval $m_a^{*}\leq m_a<\infty$ or $m_a=l_a$, then there exists a unique triangle
  with prescribed $l_a$, $l_b$ and $m_a$.
  \item Let the median $m_b$ be given. If its length belongs to the interval $l_b<m_b<m_b^{*}$,
  then there exist two triangles with prescribed $l_a$, $l_b$ and $m_b$. If $m_b$ belongs to
  the interval $m_b^{*}\leq m_b<\infty$ or $m_b=l_b$, then there exists a unique triangle
  with prescribed $l_a$, $l_b$ and $m_b$.
  \item Let the median $m_c$ be given. Since $m_c$ varies continuously from $m_c=m_c^*$ to infinity, then
  for any $m_c$ from the interval $m_c^*<m_c<\infty$ there exists a unique triangle
  with prescribed $l_a$, $l_b$ and $m_c$.
  \item Let the bisector $l_c$ be given. Since $l_c$ varies continuously from zero to infinity, then
  for any $l_c$ from the interval $0<l_c<\infty$ there exists a unique triangle
  with prescribed $l_a$, $l_b$ and $l_c$.
\end{enumerate}

If the conditions formulated in each item for the third element are not fulfilled, then the
corresponding triangle does not exist.

\section*{Acknowledgements}

I am grateful to Prof.~A.~M.~Kamchatnov for help with translation of this paper.
I thank Dr.~P.~Y.~Georgievsky and Dr.~T.~A.~Zhuravskaya for help with editing of the text
and preparation of figures.


\begin{thebibliography}{99}
\bibitem{mp-1994} P. Mironescu and L. Panaitopol, ``The existence of a triangle
with prescribed angle bisector lengths'',
Amer. Math. Monthly, vol.~101 (1994) 58-60.
\bibitem{za-2003} A. Zhukov and N. Akulich, ``Is a triangle determined uniquely?''
Kvant, N1 (2003) 29-31 (in Russian).
\bibitem{os-2016}	S. Osinkin, ``On the existence of a triangle with prescribed bisector lengths'',
Forum Geometricorum, vol.~16 (2016) 399-405.
\end{thebibliography}
\end{document}